\title{CR Embeddings, Chains, and \\the Fefferman Bundle}
\author{
        Andr\'e Minor\\
}
\date{January, 2013}
\documentclass[11pt]{article}
\usepackage{color,amsmath,amssymb}

\newtheorem{lem}{Lemma}
\newtheorem{thm}{Theorem}
\numberwithin{equation}{section}
\numberwithin{lem}{section}
\numberwithin{thm}{section}
\begin{document}
\maketitle
\nocite{CM}	
\nocite{BER99}
\nocite{EHZ}
\section{Introduction}
On a Riemannian manifold, the Levi-Civita connection induces a family of curves called geodesics defined by a certain second order partial differential equation. On a Levi-nondegenerate CR hypersurface, chains are a special type of CR invariant curve defined by a second order partial differential equation determined by the CR structure. Chains determine the CR structure of a Levi-nondegenerate CR hypersurface in the sense of the result by J. Cheng in 1988 \cite{Cheng} which says that any (local) diffeomorphism between two Levi-nondegenerate CR manifolds of hypersurface type which preserves chains must in fact be a CR (or conjugate CR) diffeomorphism. Chains are introduced in the celebrated paper of S. S. Chern and J. K. Moser \cite{CM}. A curve in a strictly pseudo-convex CR manifold M is a chain if, with respect to the forms $\{ \omega,\omega^{\alpha},\omega^{\bar{\beta}},\phi,\phi_{\alpha}^{\:\:\beta},\phi^{\alpha},\phi^{\bar{\beta}},\psi \}$ constructed in \cite{CM} (see section \ref{sec:chains}), the following equations holds along the curve:
\[\omega^{\alpha}=\phi^{\alpha}=0\]
Suppose $f$ : $ M $ $\hookrightarrow$ $\hat{M}$ is a local smooth embedding of a CR manifold M into a CR manifold $\hat{M}$ of strictly larger dimension,	 with both $M$ and $\hat{M}$ strictly pseudo-convex and of hypersurface type. In the Riemannian setting, the vanishing of the second fundamental form of a smooth embedding is equivalent to the preservation of geodesics. In the CR setting this may not generally be the case with chains, however we will see below that, in the case that the target is spherical (ie $\hat{M}$ is locally CR equivalent to the sphere $\mathbb{S}$), a CR embedding $f:M\rightarrow \mathbb{S}$ sends chains to chains if and only if its CR second fundamental form vanishes. Moreover, in this paper we will discover another geometric property of a CR embedding that is equivalent to the preservations of chains, in the case that the target is locally spherical. If the target is not spherical, we will find necessary and sufficient conditions for a CR embedding to send chains to chains. \\ 
In his 1976 paper \cite{Fef}, C. L. Fefferman constructed a Lorentz metric on the trivial circle bundle over the boundary of a strictly pseudoconvex domain in $\mathbb{C}^n$ that carried information about the CR structure of the boundary. His method involved using an approximation of a solution to the Monge-Amp$\grave{\text{e}}$re equations. In particular he proves that chains are the projections of light rays (i.e. null geodesics) on the circle bundle to the boundary of the domain. In 1977 D. Burns, Jr., K. Diederich, and S. Shnider \cite{BDS77} and S. M. Webster \cite{web77km} introduced intrinsic constructions of Fefferman's Lorentz metric, thereby generalizing the Fefferman bundle to abstract CR manifolds with nondegenerate Levi Form (not necessarily positive definite). The conformal class of the Fefferman metric is a CR invariant which we will denote by [h], where h is a metric representing the conformal class. Denote by (C $\rightarrow$ M, [h]) the Fefferman bundle. Now we can ask the following question:
\begin{quote} Under what conditions may a local CR embedding $f$ : $ M $ $\hookrightarrow$ $\hat{M}$ from a CR manifold M into a CR manifold $\hat{M}$, with both $M$ and $\hat{M}$ strictly pseudo-convex and of hypersurface type, be lifted to a conformal isometric embedding of $C$ into $\hat{C}$?
\end{quote}	
It is easy to see that any CR diffeomorphism between two equidimensional CR hypersurfaces locally lifts to a conformal isometry of the Fefferman bundles. If the target has strictly larger dimension than the source, we will show that a necessary and sufficient condition for such an isometric lift to exist will be a certain relationship between the pseudo-conformal curvature tensor of the ambient space $\hat{M}$ and the CR second fundamental form of the embedding $f$ (see Theorem \ref{techthm}). This condition will be satisfied when the embedding sends chains to chains. In addition, if the embedding sends chains to chains we will conclude that the pseudo-Riemannian second fundamental form of the lift must vanish. 
\begin{thm}
Suppose $M$ and $\hat{M}$ are strictly pseudoconvex CR hypersurfaces and $f : M \hookrightarrow \hat{M}$ is a smooth CR embedding. Then $f$ sends chains on $M$ to chains on $\hat{M}$ if and only if there is a lift of $f$ to a conformal isometric embedding of the associated Fefferman bundles with vanishing pseudo-Riemannian second fundamental form.
\label{gentargetchainthm}
\end{thm}
The celebrated theorem of Nash says that any Riemannian manifold can be isometrically embedded into Euclidean space of sufficiently large dimension. Since spheres are the locally unique strictly pseudoconvex CR hypersurfaces with vanishing pseudo-conformal curvature, they play the role of Euclidean space in CR geometry. Existence and rigidity of CR embeddings into spheres have been heavily studied but still leave many open questions. 
\begin{thm}
Suppose $f : M \hookrightarrow \mathbb{S}^{2\hat{n}+1}$ is a local smooth CR embedding of a strictly pseudoconvex smooth hypersurface $M \subset \mathbb{C}^{n+1}$ into the sphere $\mathbb{S}^{2\hat{n}+1}\subset\mathbb{C}^{\hat{n}+1}$. Let $C$ and $\hat{C}$ be the Fefferman bundles associated to $M$ and $\mathbb{S}^{2\hat{n}+1}$ respectively. The following conditions are equivalent.
\begin{enumerate}
\item $f$ sends chains on $M$ to chains on $\hat{M}$.
\vspace{.1cm}
\item There is a lift of $f$ to a conformal isometric embedding of $C$ into $\hat{C}$.
\vspace{.1cm}
\item The CR second fundamental form of f vanishes. 
\vspace{.1cm}
\item There exists a local CR diffeomorphism $\phi$ from the sphere $\mathbb{S}^{2n+1}$ to $M$ and an automorphism of the target sphere $A \in Aut(\mathbb{S}^{2\hat{n}+1})$ such that the composition $A \circ f \circ \phi$ : $ \mathbb{S}^{2n+1} \rightarrow \mathbb{S}^{2\hat{n}+1}$ is the linear embedding.
\end{enumerate}
\label{fefthmsphere}
\end{thm}
\section{Adapted Coframes}
\label{adaptedframes}
Let $M\subset\mathbb{C}^{n+1}$ be a strictly pseudoconvex CR manifold of hypersurface type. There is a subbundle $H$ of the tangent space of $M$, called the complex tangent space of $M$. A choice of a nonvanishing real 1-form $\theta \in H^\bot$ which annihilates $H$ is called a contact form on $M$. Fixing such a choice of $\theta$, the pair $ ( M, \theta ) $ is called a choice of pseudohermitian structure. We say $ \{ \theta$,  $\theta^\alpha$, $\theta^{\bar{\beta}} \}$ is an admissible coframe for $M$ if 
\begin{equation}
d\theta = ig_{\alpha \bar{\beta}} \theta^\alpha \wedge \theta^{\bar{\beta}}
\label{admissible}
\end{equation}
where $g_{\alpha \bar{\beta}}$ is a positive definite constant hermitian matrix and $1 \leq \alpha \leq n$ (notice that dim($M$) = $2n+1$). After a linear change of $\{\theta^\alpha\}$ preserving (\ref{admissible}), we may assume that $g_{\alpha \bar{\beta}}$ is the identity matrix. \\
In \cite{Web78} Webster shows that associated to a pseudohermitian structure on M is a family of connection 1-forms $\omega_\alpha^{\: \: \:\beta}$ which are uniquely determined by the choice of coframe \{$\theta$, $\theta^\alpha$, $\theta^{\bar{\beta}}$\} and the conditions:
\begin{equation}
d\theta^\beta = \theta^\alpha \wedge \omega_\alpha^{\: \: \:\beta} + \theta \wedge \tau^\beta 
\label{dthetaalpha}
\end{equation}
\begin{equation}
dg_{\alpha \bar{\beta}} = \omega_{\alpha \bar{\beta}} + \omega_{\bar{\beta} \alpha } \end{equation}
with $\tau^\beta = A_{\: \: \: \: \bar{\mu}}^{\beta} \theta^{\bar{\mu}},\: \: \: A^{\alpha \beta} = A^{\beta\alpha} $, and here we are using the summation convention with the matrix $(g_{\alpha \bar{\beta}})$ to raise and lower indices, (e.g. $\omega_{\alpha \bar{\beta}} = \omega_{\alpha}^{\:\:\: \gamma} g_{\gamma \bar{\beta}}$). The analogous forms on $\hat{M}$ will be given a hat. \\
Covariant differentiation is given by contracting with the connection 1-forms. For example, the covariant derivative of $A_{\: \: \: \: \bar{\alpha}}^{\beta}$ is given by:
\begin{equation}
\nabla A_{\: \: \: \bar{\alpha}}^{\beta}= dA_{\: \: \: \bar{\alpha}}^{\beta}+A_{\: \: \: \bar{\alpha}}^{\mu}\omega_{\mu}^{\;\;\;\beta}-A_{\: \: \: \bar{\mu}}^{\beta}\omega_{\bar{\alpha}}^{\;\;\;\bar{\mu}}= A_{\: \: \: \bar{\alpha};0}^{\beta}\theta+A_{\: \: \: \bar{\alpha};\mu}^{\beta}\theta^\mu+A_{\: \: \: \bar{\alpha};\bar{\mu}}^{\beta}\theta^{\bar{\mu}}
\end{equation}
Webster showed that these forms determine a unique connection on $H$. He then showed that the pseudohermitian curvature $R_{\alpha\: \: \: \: \: \mu \bar{\nu}}^{\: \: \: \beta}$ satisfies
\begin{equation}
d\omega_\alpha^{\:\:\:\beta} - \omega_\alpha^{\:\:\:\gamma} \wedge \omega_\gamma^{\:\:\:\beta} = R_{\alpha\: \: \: \: \: \mu \bar{\nu}}^{\: \: \: \beta} \theta^\mu \wedge \theta^{\bar{\nu}} + W_{\alpha\: \: \: \: \: \mu }^{\: \: \: \beta} \theta^\mu \wedge \theta - W_{\: \: \: \alpha \bar{\nu} }^{ \beta} \theta^{\bar{\nu}} \wedge \theta +i\theta_\alpha\wedge \tau^\beta - \tau_\alpha\wedge\theta^\beta
\label{defcurv}
\end{equation}
Given the admissible coframe \{$\theta$, $\theta^\alpha$, $\theta^{\bar{\beta}}$\} for $M$ as above it is shown in \cite{EHZ} that locally there is a choice of contact form $\hat{\theta}$ on $\hat{M}$ and an admissible coframe $\{ \hat{\theta}, \hat{\theta}^A ,\hat{\theta}^{\bar{A}}  \}$ for $\hat{M}$, where $1 \leq A \leq \hat{n}$, such that when pulled back to $M$ via the embedding $f$ : $ M $ $\hookrightarrow$ $\hat{M}$ we have
\begin{equation}
\begin{array}{l l}
\hat{\theta} = \theta \\
\hat{\theta}^\alpha = \theta^\alpha \\
\hat{\theta}^a = 0
\end{array}
\end{equation}
where we use lower case Roman letters to denote the normal direction, $n+1 \leq a \leq \hat{n}$. Such a pair of coframes is known as an adapted pair of coframes. \\
It is shown in \cite{EHZ} that on $M$ we have 
\begin{equation}
\begin{array}{cc}
\hat{\omega}_\alpha^{\: \: \:\beta} = \omega_\alpha^{\: \: \:\beta} &  \vspace{.3cm}\\
\hat{\tau}^\alpha = \tau^\alpha &\tau^b = 0  \vspace{.3cm}\\
\omega_\alpha^{\: \: \: b} = \omega_{\alpha \: \: \beta }^{\: \: \: b} \theta^\beta & \omega_{\alpha \: \: \beta }^{\: \: \: b} = \omega_{\beta \: \: \alpha }^{\: \: \: b} 
\end{array}
\end{equation}
and that the CR second fundamental form $\Pi$ is then given by:
\begin{equation}
\Pi(L_\alpha, L_\beta ) = \omega_{\alpha \: \: \beta }^{\: \: \: b} L_b
\label{Picoord}
\end{equation}
\newpage
\section{The Fefferman Bundle}
We now fix a choice of pseudohermitian structure $\theta$ and an admissible coframe $ \{ \theta$,  $\theta^\alpha$, $\theta^{\bar{\alpha}} \}$ on $M$. Let $C=M \times S^1$ and locally define the 1-form $\sigma$ by: 
\begin{equation}
\sigma = \frac{1}{n+2} (dt + i\omega_\alpha^{\: \: \alpha} - \frac{1}{2(n+1)}R\theta - \frac{i}{2} g^{\alpha \bar{\beta} } dg_{\alpha \bar{\beta}} )  
\end{equation}
where the variable t parameterizes the $S^1$ coordinate over M. Then we define the metric $h$ on C by:
\begin{equation}
h = \theta^\alpha \cdot \theta_\alpha + 2\theta \cdot \sigma 
\label{fefmet}
\end{equation}
Here the '$\cdot$' means symmetric product. It is shown in \cite{Lee86} that both the form $\sigma$ and the conformal class of the Lorentz metric $h$ are independent of the choice of admissible coframe, are globally defined on C, and that $h$ corresponds to the Fefferman metric previously developed by Fefferman (among others \cite{Fef}, \cite{BDS77}).
In [7] it is shown that the projections of null geodesic from C onto M are chains (excluding the fibers of C which project to points) and all chains are given as the projection of a null geodesic in C.
\begin{lem} With respect to this coframe we have:
\begin{equation}
R_{\alpha \bar{\beta}} - \frac{1}{2(n+1)}Rg_{\alpha \bar{\beta}} = 0 \: \: \: \iff \: \: \: i\omega_\alpha^{\: \: \: \alpha} - \frac{1}{2(n+1)}R\theta \: \: \: \text{is closed}
\end{equation}
\label{curveq}
\end{lem}
Note: The condition on the left hand side is equivalent to $R_{\alpha \bar{\beta}}=0$. 
The proof of this lemma appears in \cite{Lee88} and will be given here out of interest. It requires the following useful lemma involving the complex tangent space $H$ of $M$. 
\begin{lem}
If $\xi$ is a closed 2-form on $M$ such that $\xi | H = 0$, then $\xi = 0$.
\label{CRlem}
\end{lem}
Given a choice of contact form $\theta$ we have $H^{\perp} =\langle \theta \rangle $. The condition $\xi | H = 0$ implies that $\xi= \eta \wedge \theta$ for some 1-form $\eta$. The assumption that $\xi$ is closed then gives $0=d\eta \wedge \theta + \eta \wedge d\theta$. Restricting to $H$ we then see $\eta \wedge d\theta |_{H} = 0$. Since the Levi form is nondegenerate we then have $\eta|_{H}=0$ (ie $\eta = 0$ mod$\theta$). Thus $\xi = \eta \wedge \theta = 0$.\\ QED 
\vspace{.5cm}
\\ 
Now to prove Lemma \ref{curveq} we observe from equation (\ref{defcurv}) and by the symmetry properties of the connection forms that we have: 
\begin{equation}
d\omega_\alpha^{\: \: \: \alpha} = d\omega_\alpha^{\: \: \: \alpha} - \omega_\alpha^{\: \: \: \gamma} \wedge \omega_\gamma^{\: \: \: \alpha} = R_{\mu \bar{\nu}}\theta^\mu \wedge \theta^{\bar{\nu}} \text{ mod } \theta
\label{domega}
\end{equation}
Thus we see $R_{\alpha \bar{\beta}} - \frac{1}{2(n+1)}Rg_{\alpha \bar{\beta}} = 0$ if and only if 
\begin{equation}
d\omega_\alpha^{\: \: \: \alpha} = -\frac{i}{2(n+1)} R d\theta = -\frac{i}{2(n+1)} d(R\theta) \text{ mod } \theta 
\end{equation}
which, by lemma \ref{CRlem}, holds if and only if
\begin{equation}
d\omega_\alpha^{\: \: \: \alpha} +\frac{i}{2(n+1)} d(R\theta) = 0 
\end{equation}
QED
\section{Chains}
\label{sec:chains}
Let $ \{ \theta$,  $\theta^\alpha$, $\theta^{\bar{\alpha}} \}$ be an admissible coframe on a smooth strictly pseudoconvex CR hypersurface $M \subset \mathbb{C}^{n+1}$. As in \cite{CM}, we let $E$ be the ray bundle of oriented contact forms on $M$. $\{u\theta \; |\; u \geq 0 \text{ is smooth on } M \}$ is the space of sections of this bundle. The form $\omega = u\theta$ is intrinsically defined on $E$. By equation (\ref{admissible}), we have
\begin{equation}
\begin{array}{cc}
d\omega &= iug_{\alpha \bar{\beta}}\theta^\alpha \wedge\theta^{\bar{\beta}} + du \wedge \theta\\  
& = ig_{\alpha \bar{\beta}}\omega^\alpha \wedge\omega^{\bar{\beta}} + \omega\wedge \phi
\end{array}
\label{cmadmissible}
\end{equation}
where $\omega^\alpha = \sqrt{u}\theta^\alpha$ and $\phi=-\frac{du}{u}$. Define $G_1$ to be the space of linear transformations on $(\omega,\omega^\alpha,\omega^{\bar{\beta}},\phi)$ that preserve $\omega$ and equation (\ref{cmadmissible}). Then we let $Y$ be the principal $G_1$ bundle of coframes $(\omega,\omega^\alpha,\omega^{\bar{\beta}},\phi)$ on $E$, $Y\rightarrow E \rightarrow M$. By (\ref{admissible}), $g_{\alpha \bar{\beta}}$ is a positive definite constant hermitian matrix. In $Y$ we have a family of linear differential forms 
\begin{equation}
\{ \omega,\omega^{\alpha},\omega^{\bar{\beta}},\phi,\phi_{\alpha}^{\:\:\beta},\phi^{\alpha},\phi^{\bar{\beta}},\psi \}
\end{equation}
determined by the following structure equations from \cite{CM}
\begin{equation}
\begin{array}{c}
\phi_{\alpha \bar{\beta}}+\phi_{\bar{\beta}\alpha} -g_{\alpha \bar{\beta}}\phi = 0 \vspace{.2cm}\\ 
d\omega=ig_{\alpha \bar{\beta}}\omega^\alpha\wedge\omega^{\bar{\beta}}+\omega\wedge\phi \vspace{.2cm}\\
d\omega^\alpha=\omega^\beta\wedge\phi_\beta^{\;\;\alpha}+\omega\wedge\phi^\alpha \vspace{.2cm}\\
d\phi = i \omega_{\bar{\beta}}\wedge\phi^{\bar{\beta}}+i\phi_{\bar{\beta}}\wedge \omega^{\bar{\beta}} +\omega \wedge \psi \vspace{.2cm}	\\
d\phi_\beta^{\;\; \alpha}=\phi_\beta^{\;\; \sigma}\wedge \phi_\sigma^{\;\; \alpha}+i\omega_\beta\wedge\phi^\alpha-i\phi_\beta\wedge\omega^\alpha-i\delta_\beta^{\;\;\alpha}(\phi_\sigma\wedge\omega^\sigma)-\frac{1}{2}\delta_\beta^{\;\;\alpha}\psi\wedge\omega+\Phi_\beta^{\;\;\alpha} \vspace{.2cm}	\\
d\phi^\alpha=\phi\wedge\phi^\alpha+\phi^\beta\wedge\phi_\beta^{\;\;\alpha}-\frac{1}{2}\psi\wedge\omega^\alpha+\Phi^\alpha \vspace{.2cm}	\\
d\psi=\phi\wedge\psi+2i\phi^\beta\wedge\phi_\beta+\Psi
\end{array}
\label{structureeqns}
\end{equation}
The curvature forms are given by:
\[
\begin{array}{c}
\Phi_\beta^{\;\;\alpha}=S_{\beta \; \rho \bar{\sigma}}^{\;\;\alpha}\omega^\rho\wedge\omega^{\bar{\sigma}}+V_{\beta \; \rho}^{\;\alpha}\omega^\rho\wedge\omega-V_{\;\;\beta  \bar{\sigma}}^{\alpha}\omega^{\bar{\sigma}}\wedge\omega \vspace{.2cm}	\\
\Phi^\alpha=-V_{\rho \; \sigma}^{\;\alpha}\omega^\rho\wedge\omega^\sigma+V_{\;\;\rho\bar{\sigma}}^{\alpha}\omega^\rho\wedge\omega^{\bar{\sigma}}+P_\rho^{\;\;\alpha}\omega^\rho\wedge\omega+Q_{\bar{\sigma}}^{\;\;\alpha}\omega^{\bar{\sigma}}\wedge\omega \vspace{.2cm}	\\
\Psi = iQ_{\rho \sigma}\omega^\rho\wedge\omega^\sigma-iQ_{\bar{\rho}\bar{\sigma}}\omega^{\bar{\rho}}\wedge\omega^{\bar{\sigma}}-i(P_{\rho\bar{\sigma}}+P_{\bar{\sigma}\rho})\omega^\rho\wedge\omega^{\bar{\sigma}}+(R_\rho\omega^\rho+R_{\bar{\sigma}}\omega^{\bar{\sigma}})\wedge\omega
\end{array}
\]
The forms $\phi^\alpha, \phi^{\bar{\beta}},\phi_\alpha^{\;\;\beta}, \psi$ are uniquely determined by $(\omega,\omega^\alpha,\omega^{\bar{\beta}},\phi)$ by requiring the pseudo-conformal curvature tensors of $M$ satisfy:
\[
\begin{array}{c}
S_{\alpha\bar{\beta}\rho\bar{\sigma}}=S_{\rho\bar{\beta}\alpha\bar{\sigma}}=S_{\alpha\bar{\sigma}\rho\bar{\beta}} \vspace{.2cm}	\\
S_{\alpha\bar{\beta}\rho\bar{\sigma}}=S_{\bar{\beta}\alpha\bar{\sigma}\rho} (=\overline{S_{\beta\bar{\alpha}\sigma\bar{\rho}}}) \vspace{.2cm}	\\
S_{\rho \; \alpha \bar{\beta}}^{\;\;\rho}=V_{\beta \; \rho}^{\;\;\rho}=P_{\rho}^{\;\;\rho}=0
\end{array}
\]
A choice of an admissible coframe $\{ \theta$, $\theta^\alpha$, $\theta^{\bar{\alpha}} \}$ (where $\theta\in H^{\perp}, \theta \neq 0$ is a real 1-form), induces a section $\{\omega,\omega^\alpha,\omega^{\bar{\beta}},\phi\}$ of $Y$ which we use to pull forms on $Y$ back to $M$. The pullback of the forms $\{\omega,\omega^\alpha,\omega^{\bar{\beta}},\phi\}$ is $\{\theta, \theta^\alpha,\theta^{\bar{\beta}},0\}$. We get an induced family of 1-forms on $M$ \{$\phi_\beta^{\: \: \: \alpha}$,
$\phi^{ \alpha}$, $\phi^{\bar{\alpha}}$, $\psi \}$ which are the pullbacks of the associated forms on $Y$.
\\
A curve $\gamma$ which is transverse to the CR tangent space (ie $\theta(\dot{\gamma})\neq0$ for any $0\neq\theta \in H^{\bot}$) is a chain if, when we choose our admissible coframe $(\theta,\theta^\alpha,\theta^{\bar{\beta}})$ so that along $\gamma$ 
\begin{equation}
\theta^\alpha = 0,
\label{alpha0}
\end{equation}
we also have,
\begin{equation}
\phi^\alpha = 0.
\label{phiiszero}
\end{equation}
It is not hard to show that this definition is independent of the choice \{$\theta^\alpha$, $\theta^{\bar{\beta}}\}$ satisfying (\ref{alpha0}) (so that $ \{ \theta$,  $\theta^\alpha$, $\theta^{\bar{\beta}} \}$ is an admissible coframe).
\\ 
The above definition can be complicated to work with, because for a given curve one must choose a coframe associated to that curve satisfying (\ref{alpha0}). When dealing with more than one chain, a geometric comparison of the two curves is then very difficult. In the proof of \cite{Cheng}, another formulation of the definition of a chain is given which may be used with any admissible coframe. 
\begin{lem}
Suppose $M$ is a strictly pseudoconvex CR manifold of hypersurface type. A curve $\gamma$ in $M$ is a chain if and only if for any admissible coframe $\{ \theta$, $\theta^\alpha$, $\theta^{\bar{\alpha}} \}$ the following equation has a solution $a(t)=(a^1,...,a^n)$ along $\gamma$:
\begin{equation}
\left\{ \begin{array}{c c}
\theta^\alpha = 2a^\alpha \theta
\\
da^\alpha = 4ia^\alpha |a|^2 \theta -  a^\beta\phi_\beta^{\: \: \:\alpha}-\frac{1}{2} \phi^\alpha
\end{array} \right.
\label{13}
\end{equation}
\label{psxaltdef}
\end{lem}
Proof of lemma \ref{psxaltdef}: Fix a chain $\gamma \subset M$. Fix a contact form $\theta$ on $M$ and an admissible coframe $\{ \theta, \theta^\alpha \}$. After a change of variables, the admissibility condition means
\[
d\theta = i g_{\alpha \bar{\beta}} \theta^\alpha \wedge \theta^{\bar{\beta}}
\]
were the $n\times n$ matrix $g_{\alpha \bar{\beta}}$ is the identity matrix representing the Levi Form. The coframe constructed in \cite{CM} may be written in the form of a Maurer-Cartan form as:
\begin{equation}
\pi = \left(  \begin{array}{c c c}
\pi_0^{\: \: 0} & \theta^\alpha & 2 \theta \\
-i\phi_\alpha & (\phi_\alpha^{\: \: \:\beta} + \pi_0^{\: \: 0}
\delta_\alpha^{\: \:
\:\beta}) & 2i \theta_\alpha  \\
-\frac{1}{4} \psi & \frac{1}{2} \phi^\alpha &-\bar{\pi}_0^{\: \: 0}
\end{array} \right)
\end{equation}
where $\pi_0^{\: \: 0} = -\frac{1}{n+2} (\phi_\alpha^{\: \: \:\alpha}+\phi)$. This family of forms is unique up to transformation by an element of the G-structure, $h\in G$. The transformation is then given on $M$ by:
\begin{equation}
\tilde{\pi} = dh h^{-1} + h \pi h^{-1} 
\label{trnsflaw}
\end{equation}
We are free to choose a coframe $\theta$, $\theta^\alpha$ so that along $\gamma$ we have $\theta = 1$. It is shown in \cite{BDS77} that the element $h \in G$ may be chosen to be of the form: 
\begin{equation}
h(t) = \left( \begin{array}{c c c}
1 & 0 & 0 \\
-2ia^* & I & 0 \\
-i|a|^2 & a & 1
\end{array} \right) \: \: \:, a(t) \in \mathbb{C}^n
\end{equation}
$\gamma$ is a chain if and only if there is some transformation h so that $\tilde{\theta}^\alpha=0=\tilde{\phi}^\alpha$. One may now isolate the $\tilde{\theta}^\alpha$, $\tilde{\phi}^\alpha$ terms from equation (\ref{trnsflaw}) and set them equal to 0. Using:
\begin{equation}
h^{-1}(t) = \left( \begin{array}{c c c}
1 & 0 & 0 \\
2ia^* & I & 0 \\
-i|a|^2 & -a & 1
\end{array} \right) \: \: \:, a(t) \in \mathbb{C}^n
\end{equation}
and letting $I,J \in {0,...,n+1}$ we have:
\begin{equation}
\begin{array}{rl}
2\tilde{\theta} &=  \tilde{\pi}_0^{\:\:\: n+1} = (dh)_{0}^{\:\:\:I} (h^{-1})_{I}^{\:\:\: n+1} + h_{0}^{\:\:\:I} \pi_I^{\:\:J} (h^{-1})_J^{\:\:\:n+1} 
\vspace{.2cm}
\\ & = 2 \theta
\vspace{.8cm}
\\
\tilde{\theta}^\alpha &= \tilde{\pi}_0^{\:\:\: \alpha}  =  (dh)_{0}^{\:\:\:I} (h^{-1})_{I}^{\:\:\: \alpha} + h_{0}^{\:\:\:I} \pi_I^{\:\:J} (h^{-1})_J^{\:\:\:\alpha} 
\vspace{.2cm}
\\ & = (\phi_0^{\:\:\:0}, \theta^\beta, 2\theta) \cdot (0, \delta_\beta^{\:\:\:\alpha}, -a^\alpha) 
\vspace{.2cm}
\\ & = \theta^\alpha - 2a^\alpha \theta
\vspace{.8cm}\\
\frac{1}{2} \tilde{\phi}^{\alpha} & = \tilde{\pi}_{n+1}^{\:\:\: \alpha} 
\vspace{.2cm}\\ 
& = (dh)_{n+1}^{\:\:\:I} (h^{-1})_{I}^{\:\:\: \alpha} + h_{n+1}^{\:\:\:I} \pi_I^{\:\:J} (h^{-1})_J^{\:\:\:\alpha} 
\vspace{.2cm}\\
& = da^\alpha \\
&\;\; + (-i|a|^2, a^\beta, 1) \cdot (\theta^\alpha -2a^\alpha, \phi_\beta^{\:\:\:\alpha} + \pi_0^{\:\:0} \delta_\beta^{\:\:\:\alpha} - 2ia^\alpha \theta_\beta, \frac{1}{2} \phi^\alpha + a^\alpha  \bar{\pi}_0^{\:\:0} ) 
\vspace{.2cm}\\ 
& = da^\alpha + a^\beta \phi_\beta^{\:\:\:\alpha} + a^\beta \pi_0^{\:\:0} \delta_{\beta}^{\:\:\:\alpha} - 2ia^\alpha a^\beta \theta_\beta + \frac{1}{2}\phi^\alpha +a^\alpha \bar{\pi}_0^{\:\:0}
\vspace{.2cm}\\ 
& = da^\alpha- 2ia^\alpha a^\beta \theta_\beta + a^\beta \phi_\beta^{\:\:\:\alpha} + a^\alpha (\pi_0^{\:\:0} +\bar{\pi}_0^{\:\:0} ) + \frac{1}{2}\phi^\alpha 
\end{array}
\end{equation}
\vspace{.5cm}
\\
Note: When we choose an admissible coframe we have $\phi = 0$ on $M$. Since $\pi_0^{\: \: 0} = -\frac{1}{n+2} (\phi_\alpha^{\: \: \:\alpha}+\phi)$ and by the anti-symmetry property of the forms $\phi_\alpha^{\:\:\:\beta}$ (\ref{structureeqns}), we conclude $\pi_0^{\:\:0}+ \bar{\pi}_0^{\: \: 0} = 0$. Setting $\tilde{\theta}^\alpha = \tilde{\phi}^\alpha = 0$ results in the desired differential equation on $a(t)$. QED
\newpage
\section{An Isometric Lift}
\label{Lift}
\begin{thm}
The map $f : M \hookrightarrow \hat{M}$ may be locally lifted to a conformal isometric embedding of the Fefferman bundles if and only if the following equation holds with respect to any adapted coframe.
\begin{equation}
\hat{S}_{a \: \: \: \: \alpha \bar{\beta}}^{\: \: \: a} + \omega_{\mu \:
\: \: \alpha}^{\: \: \: a}\omega_{\: \: \: a \bar{\beta} }^\mu =0
\label{techeqn}
\end{equation}
\label{techthm}
\end{thm}
Equation (\ref{techeqn}) is a condition on part of the CR conformal curvature tensor $\hat{S}$ of $\hat{M}$ and the second fundamental form of the map $f$. Thus it is a condition that is independent of the choice of pseudohermitian structure. \\
Let $\{ \hat{\theta}, \hat{\theta}^A , \hat{\theta}^A \}$ be an admissible coframe on $\hat{M}$ adapted to the admissible coframe $\{ \theta, \theta^\alpha , \theta^{\alpha} \}$ on $M$. Recall by equation (\ref{fefmet}), the conformal classes of the Fefferman metrics on $C$ and $\hat{C}$ are represented by the metrics 
\begin{equation}
\begin{array}{c}
h = \theta^\alpha \cdot \theta_\alpha + 2\theta \cdot \sigma
\vspace{.2cm} \\ 
\hat{h} = \hat{\theta}^A \cdot \hat{\theta}_A + 2\hat{\theta}\cdot \hat{\sigma}
\end{array}
\end{equation}
respectively. The form $\sigma$ is given by
\begin{equation}
\begin{array}{c}
\sigma = \frac{1}{n+2} (dt + i\omega_\alpha^{\: \: \alpha} - \frac{1}{2(n+1)}R\theta - \frac{i}{2} g^{\alpha \bar{\beta} } dg_{\alpha \bar{\beta}} ) 
\vspace{.2cm}
\\
\hat{\sigma} = \frac{1}{\hat{n}+2} (ds + i\hat{\omega}_A^{\: \: A} - \frac{1}{2(\hat{n}+1)}\hat{R}\hat{\theta} - \frac{i}{2} \hat{g}^{A \bar{B} } d\hat{g}_{A \bar{B}} ) 
\end{array}
\label{sighatsig}
\end{equation}
We will take the variable s for the $S^1$ coordinate in $\hat{C}$. With respect to this adapted coframe, since $f^{*}( \hat{\theta}^A \cdot \hat{\theta}_A) = \theta^\alpha \cdot \theta_\alpha$ and $f^{*} \hat{\theta}=\theta$, any lift $F:C\hookrightarrow \hat{C}$ of the CR embedding $f:M\hookrightarrow \hat{M}$ is a conformal isometry if and only if it is actually an isometry and $f^* \hat{\sigma}=\sigma$.
\\
By \cite{EHZ} a relationship between Webster's connection 1-forms and the pullbacks of the Chern-Moser forms is given by:
\begin{equation}
\phi_{\beta}^{\: \: \: \alpha} = \omega_{\beta}^{\: \: \: \alpha} + D_{\beta}^{\: \: \: \alpha}\theta, \: \: \: \: \: \: \phi^{\: \alpha} = \tau^\alpha + D_{\mu}^{\: \: \: \alpha}\theta^\mu + E^\alpha\theta 
\label{CMWeb}
\end{equation}
where
\begin{equation}
\begin{array}{l l}
D_{\alpha \bar{\beta}} = \frac{i}{n+2} R_{\alpha \bar{\beta}}-\frac{i}{2(n+1)(n+2)}R g_{\alpha \bar{\beta}}
\\ \\
E^\alpha = \frac{2i}{2n+1}(A_{\: \: \: \: \: ; \mu}^{\alpha \mu} - D_{\:\: \: \: \: ; \bar{\nu}}^{\bar{\nu} \alpha})
\end{array}
\label{DandE}
\end{equation}
The pullbacks of the associated forms on $\hat{M}$ are related by:
\begin{equation}
\hat{\phi}_{\beta}^{\: \: \: \alpha} = \phi_{\beta}^{\: \: \: \alpha}+C_{\beta}^{\: \: \: \alpha} \theta,\: \: \: \: \: \: \hat{\phi}^{\: \alpha} = \phi^{\: \alpha} + C_{\mu}^{\:\: \: \alpha}\theta^\mu + F^\alpha \theta
\label{CMpullbacks}
\end{equation}
where
\begin{equation}
C_{\beta}^{\: \: \: \alpha} = \hat{D}_{\beta}^{\: \: \: \alpha} -D_{\beta}^{\: \: \: \alpha}, \: \: \: \: \: \: F^\alpha = \hat{E}^\alpha -
E^\alpha
\end{equation}
It is shown in \cite{EHZ} using the pseudohermitian Gauss equation that:
\begin{equation}
C_{\alpha \bar{\beta}} = \frac{i}{n+2}(\hat{S}_{a \: \: \: \: \alpha\bar{\beta}}^{\: \: \: a} + \omega_{\mu \: \: \: \alpha}^{\: \: \:
a}\omega_{\: \: \: a \bar{\beta} }^\mu - \frac{1}{2(n+1)}(\hat{S}_{a \: \:\: \mu}^{\: \: \: a \: \;\;  \mu} + \omega_{\mu \: \: \: \nu}^{\: \: \:
a}\omega_{\: \: \: a}^{\mu \: \: \: \nu})g_{\alpha \bar{\beta}})
\label{Cs}
\end{equation}
Using equation (\ref{sighatsig}) the condition $f^* \hat{\sigma}=\sigma$ then becomes:
\begin{equation}
\frac{1}{\hat{n} +2} ds = \frac{1}{n+2}dt
+\frac{1}{n+2}(i\omega_\alpha^{\: \: \: \alpha} - \frac{1}{2(n+1)}R\theta)
- \frac{1}{\hat{n}+2}(i\hat{\omega}_A^{\: \: \: A} -
\frac{1}{2(\hat{n}+1)}\hat{R}\theta)
\label{isocond}
\end{equation}
This is a differential equation in $s=s(x,t)$ (where x denotes local coordinates on M). Locally a form is closed if and only if it is exact, thus we may locally solve equation (\ref{isocond}) if and only if the term on the right hand side is closed. We observe,
\begin{equation}
\begin{array}{l l}
&\frac{1}{n+2}(i\omega_\alpha^{\: \: \: \alpha} - \frac{1}{2(n+1)}R\theta)
- \frac{1}{\hat{n}+2}(i\hat{\omega}_A^{\: \: \: A} -
\frac{1}{2(\hat{n}+1)}\hat{R}\theta) \text{ is closed}
\\ \\
\Leftrightarrow & d(\frac{1}{n+2}\omega_\alpha^{\: \: \: \alpha} -
\frac{1}{\hat{n}+2}\hat{\omega}_A^{\: \: \: A}  ) = - i(
\frac{1}{2(n+1)(n+2)}R - \frac{1}{2(\hat{n}+1)(\hat{n}+2)}\hat{R} )
d\theta \text{ mod }\theta
\\ \\
& (\text{by Lemma \ref{CRlem}}) 
\\ \\
\Leftrightarrow & (\frac{1}{n+2}R_{\alpha \bar{\beta}} -
\frac{1}{\hat{n}+2}\hat{R}_{\alpha \bar{\beta}})= (
\frac{1}{2(n+1)(n+2)}R - \frac{1}{2(\hat{n}+1)(\hat{n}+2)}\hat{R}
)g_{\alpha \bar{\beta}}
\\ \\
& (\text{by equation (\ref{domega}) and since } d\theta = ig_{\alpha \bar{\beta}} \theta^\alpha \wedge \theta^{\bar{\beta}})
\\ \\
\Leftrightarrow & C_{\alpha \bar{\beta}} = \hat{D}_{\alpha \bar{\beta}} -
D_{\alpha \bar{\beta}}=0
\end{array}
\label{iffeqn}
\end{equation}
This means that an isometric lift of $f$ exists if and only if we have $\hat{\phi}_\beta^{\: \: \: \alpha}=\phi_\beta^{\: \: \: \alpha}$ with respect to any adapted coframe. By (\ref{Cs}) this is equivalent to:
\begin{equation}
\hat{S}_{a \: \: \: \: \alpha \bar{\beta}}^{\: \: \: a} + \omega_{\mu \: \: \: \alpha}^{\: \: \: a}\omega_{\: \: \: a \bar{\beta} }^\mu =
\frac{1}{2(n+1)}(\hat{S}_{a \: \: \: \mu}^{\: \: \: a \: \: \: \mu} + \omega_{\mu \: \: \: \nu}^{\: \: \: a}\omega_{\: \: \: a}^{\mu \: \: \:
\nu})g_{\alpha \bar{\beta}}
\label{cis0ass}
\end{equation}
By contracting both sides this becomes:
\begin{equation}
\hat{S}_{a \: \: \: \mu}^{\: \: \: a\;\;\mu} + \omega_{\mu \: \: \: \nu}^{\: \: \: a}\omega_{\: \: \: a }^{\mu\;\;\;\nu} =
\frac{n}{2(n+1)}(\hat{S}_{a \: \: \: \mu}^{\: \: \: a \: \: \: \mu} + \omega_{\mu \: \: \: \nu}^{\: \: \: a}\omega_{\: \: \: a}^{\mu \: \: \:
\nu})
\end{equation}
and thus 
\begin{equation}
\hat{S}_{a \: \: \: \mu}^{\: \: \: a\;\;\mu} + \omega_{\mu \: \: \: \nu}^{\: \: \: a}\omega_{\: \: \: a }^{\mu\;\;\;\nu}=0
\end{equation}
Plugging this back into equation equation (\ref{cis0ass}), we get
\begin{equation}
\hat{S}_{a \: \: \: \: \alpha \bar{\beta}}^{\: \: \: a} + \omega_{\mu \: \: \: \alpha}^{\: \: \: a}\omega_{\: \: \: a \bar{\beta} }^\mu=0
\label{noncontracted}
\end{equation}
Clearly if (\ref{noncontracted}) holds, then so does (\ref{cis0ass}), thus the theorem is proved. 
\\
QED 
\vspace{.1in}
\\
Notice that if $D_{\alpha \bar{\beta}}=\displaystyle R_{\alpha \bar{\beta}} - \frac{1}{2(n+1)}Rg_{\alpha \bar{\beta}} = 0 $, then by contracting we get 
\[
\displaystyle 0=R - \frac{n}{2(n+1)}R=\frac{n+2}{2(n+1)}R
\]
and thus $\displaystyle R=0$. Plugging this back into $D_{\alpha \bar{\beta}}$ we see that
\begin{equation}
D_{\alpha \bar{\beta}}=R_{\alpha \bar{\beta}} - \frac{1}{2(n+1)}Rg_{\alpha \bar{\beta}} = 0  \text{ if and only if } R_{\alpha \bar{\beta}} = 0 
\end{equation}
It is clear by equation (\ref{iffeqn}) that if the embedding $f:M\rightarrow \hat{M}$ admits an adapted pair of admissible coframes so that $R_{\alpha \bar{\beta}}=0$ and $\hat{R}_{A \bar{B}}=0$, then $C_{\alpha \bar{\beta}}=0$ and thus $f$ may be lifted to an conformal isometry of the Fefferman bundles.
\section{Comparing the Second Fundamental Forms}
It is no surprise that the second fundamental forms of $f$ and a conformal isometric lift F are related in a simple way. Fix an adapted coframe as before. Later we will use this technical lemma:
\begin{lem}
\begin{equation}
\nabla \omega_{\alpha \: \: \: \gamma}^{\:\:\: a} =
-\hat{R}_{\alpha \: \: \: \gamma \bar{\nu}}^{\:\:\: a}\theta^{\bar{\nu}} \text{
mod} \theta, \: \theta^\beta
\end{equation}
\label{covlem}
\end{lem}
To prove this lemma we start with the following identity from \cite{EHZ}:
\begin{equation}
\nabla \omega_{\alpha \: \: \: \gamma}^{\:\:\: a} = d\omega_{\alpha \: \: \: \gamma}^{\:\:\: a} - \omega_{\mu \: \: \:
\gamma}^{\:\:\: a}\omega_\alpha^{\: \: \: \mu} + \omega_{\alpha \: \: \:\gamma}^{\:\:\: b}\hat{\omega}_b^{\: \: \: a}- \omega_{\alpha \: \: \:
\mu}^{\:\:\: a}\omega_\gamma^{\: \: \: \mu}
\label{covderivpi}
\end{equation}
This is the covariant derivative of $\omega_{\alpha \: \: \: \gamma}^{\:\:\: a}$ as a section of the bundle of $\mathbb{C}$-bilinear maps 
\[
T_p^{1,0}M \times T_p^{0,1}M \rightarrow T_p^{1,0}\hat{M} /T_p^{1,0}M.  
\]
See \cite{EHZ} for more details. Working mod$\theta$, pulling back to M, and using the $\hat{M}$ analog of equations (\ref{defcurv}), (\ref{Picoord}), (\ref{dthetaalpha}), (\ref{covderivpi}) we have:
\begin{align*}
-\hat{R}_{\alpha \: \: \: \mu \bar{\nu}}^{\:\:\: a}\theta^\mu \wedge \theta^{\bar{\nu}} = & \omega _{\alpha}^{\:\:\: \gamma} \wedge \omega_{\gamma}^{\:\:\: a}+ \hat{\omega} _{\alpha}^{\:\:\: b} \wedge \hat{\omega}_{b}^{\:\:\: a} - d\hat{\omega} _{\alpha}^{\:\:\: a} 
\\
= & \omega _{\alpha}^{\:\:\: \gamma} \wedge \omega _{\gamma \: \: \: \: \beta}^{\:\:\: \alpha}\theta^\beta+ \omega _{\alpha \: \: \: \beta}^{\:\:\: b}\theta^\beta \wedge \hat{\omega}_{b}^{\:\:\: a} - d(\omega_{\alpha \:\:\: \beta}^{\:\:\: a}\theta^\beta) 
\\
= & \theta^\beta \wedge (\omega _{\alpha \: \: \: \beta}^{\:\:\: b}\hat{\omega}_{b}^{\:\:\: a} - \omega _{\gamma \: \: \: \: \beta}^{\:\:\:\alpha} \omega_{\alpha}^{\:\:\: \gamma} - \omega_{\alpha \: \: \: \:\gamma}^{\:\:\: a} \omega _{\beta}^{\:\:\: \gamma}+d\omega_{\alpha \: \:\: \: \beta}^{\:\:\: a}) \\
= & \theta^\beta \wedge \nabla \omega_{\alpha \: \: \:\beta}^{\:\:\: a} \text{ mod }\theta
\end{align*}
This ends the proof of Lemma \ref{covlem}	. QED
\\
Using the coframe \{$\theta$, $\theta^\alpha$, $\theta^{\bar{\alpha}}$\} on $M$ we may take $\{\xi^I\} =$ \{$\theta$, $\sigma$, $\theta^\alpha$,  $\theta^{\bar{\alpha}}$\}, $I=-1,0,1,...,2n$, as a coframe on C and \{$T$, $X$, $L_\alpha$, $L_{\bar{\alpha}}$ \} as a dual coframe. The metric h is then given by:
\begin{equation}
h =
\left(
\begin{array}{c c c c}
0 & 1 & 0 & 0 \\
1 & 0 & 0 & 0 \\
0 & 0 & 0 & \frac{1}{2}g_{\alpha \bar{\beta}} \\
0 & 0 & \frac{1}{2}g_{\bar{\alpha} \beta  } & 0
\end{array}
\right)
\end{equation}
The connection 1-forms, $\sigma_I^{\;\;J}$ of this metric on $C$ are uniquely determined by the equations:
\begin{equation}
\begin{array}{c}
d\xi^I=\xi^J\wedge \sigma_J^{\;\;I} \\
\sigma_I^{\;\;K}h_{K J}+\sigma_J^{\;\;K}h_{K I}=dh_{I J}
\end{array}
\label{con1forms}
\end{equation}
In \cite{Lee86}, Lee finds equations for these connection 1-forms. 
\begin{lem}
The Levi-Civita connection 1-forms of h are given by:
\begin{equation}
\Omega =
\left(
\begin{array}{c c c c}
0 & 0 & i\sigma^\alpha & -i\sigma^{\bar{\alpha}} \\
0 & 0 & i\theta^\alpha & -i\theta^{\bar{\alpha}} \\
\frac{i}{2}\theta_\beta & \frac{i}{2}\sigma_\beta & \sigma_{\beta}^{\: \:
\: \alpha}  & 0  \\
-\frac{i}{2}\theta_{\bar{\beta}} & -\frac{i}{2}\sigma_{\bar{\beta}} & 0  & \sigma_{\bar{\beta}}^{\: \: \: \bar{\alpha}}  \\
\end{array}
\right)
\end{equation}
where
\[
\sigma_{\beta}^{\: \: \: \alpha} = \omega_{\beta}^{\: \: \: \alpha} + D_{\beta}^{\: \: \: \alpha} \theta + i\delta_{\beta}^{\: \: \: \alpha}
\sigma = \phi_{\beta}^{\: \: \: \alpha} + i \delta_{\beta}^{\: \: \:\alpha} \sigma
\]
\[
\sigma_\beta = i\tau_\beta -i D_{\beta \bar{\gamma}} \theta^{\bar{\gamma}}+ C_\beta \theta \text{,} \: \: \: \: \: \: C_\beta=\frac{2}{n+2}(A_{\alpha \beta ;}^{\: \: \: \: \: \: \alpha} +\frac{i}{2(n+1)} R_{\; ; \; \beta})
\]
\end{lem}
It is not difficult to show that these forms satisfy (\ref{con1forms}). The analogous procedure using the coframe \{$\hat{\theta}$, $\hat{\sigma}$, $\hat{\theta}^A$,  $\hat{\theta}^{\bar{A}}$\} on $\hat{C}$ gives:
\begin{equation}
\hat{\Omega} =
\left(
\begin{array}{c c c c}
0 & 0 & i\hat{\sigma}^A & -i\hat{\sigma}^{\bar{A}} \\
0 & 0 & i\hat{\theta}^A & -i\hat{\theta}^{\bar{A}} \\
\frac{i}{2}\hat{\theta}_B & \frac{i}{2}\hat{\sigma}_B & \hat{\sigma}_{B}^{\: \:
\: A}  & 0  \\
-\frac{i}{2}\hat{\theta}_{\bar{B}} & -\frac{i}{2}\hat{\sigma}_{\bar{B}} & 0  &
\hat{\sigma}_{\bar{B}}^{\: \: \: \bar{A}}  \\
\end{array}
\right)
\end{equation}
Since $F^* \hat{\sigma} = \sigma$ the second fundamental form is then determined by the pullback via F of the forms: 
\begin{equation}
\{\hat{\Omega}_{I}^{\: \: k} \} = \{
i\hat{\sigma}^a,\: \: -i\hat{\sigma}^{\bar{a}}, \: \: i\theta^a, \: \: 	-i\theta^{\bar{a}}, \: \: \hat{\sigma}_\alpha^{\: \: \: a}, \: \: \hat{\sigma}_{\bar{\alpha}}^{\: \: \: \bar{a}}
\}
\end{equation}
where $I$ represents tangential directions to $C$ and $k$ represents normal directions to $C$. On $C$ we have $\tau^a = 0 = \theta^a$ and thus by equation (\ref{dthetaalpha}), $\hat{A}_{\;\;\bar{\alpha}}^{a} = 0$ on $M$ and thus $\hat{A}_{\;\;\bar{\alpha} ;}^{a\;\;\;\bar{\alpha}}= 0$ on $M$. It is also clear that $\hat{D}_{\: \: \: \gamma}^a=\frac{i}{\hat{n}+2}\hat{R}_{\: \: \: \gamma}^a$. Thus the non-trivial terms in the second fundamental form of F are
\begin{equation}
\left\{
\begin{array}{l}
i\hat{\sigma}^a = -\frac{i}{\hat{n}+2} \hat{R}_{\: \: \: \gamma}^a\theta^\gamma +\frac{2}{\hat{n}+2}(\hat{A}_{\bar{b}\;\; ;}^{\;\;\;a\;\bar{b}} -\frac{i}{2(\hat{n}+1)}\hat{R}_{\: ;}^{\: \: a})\theta
\\ \\
\hat{\sigma}_\alpha^{\;\;\;a}=\hat{\omega}_\alpha^{\: \: \: a} + \frac{i}{\hat{n}+2} \hat{R}_\alpha^{\:\: \: a}\theta
\end{array}
\right.
\label{omegaterms}
\end{equation}
and their complex conjugates. We are now in a position to prove the following lemma.
\begin{lem} Let $f:M\rightarrow\hat{M}$ and $F:C\rightarrow\hat{C}$ be a conformal lift of the associated Fefferman bundles. 
\begin{enumerate}
\item If the second fundamental form of $F$ vanishes, the CR second fundamental form of $f$ vanishes. 
\item The second fundamental form of $F$ vanishes if and only if on $M$
\begin{equation}
\left\{
\begin{array}{c}
\omega_{\alpha \;\; \beta}^{\;\;a} = 0 
\vspace{.2cm}
\\
\hat{R}_\alpha^{\;\;a}=0
\vspace{.1cm}
\\
\hat{A}_{\bar{b}\;\; ;}^{\;\;a\;\;\bar{b}} -\frac{i}{2(\hat{n}+1)}\hat{R}_{\: ;}^{\: \: a} = 0
\end{array}
\right.
\end{equation}
\end{enumerate}
\label{2ndffthm}
\end{lem}
These conditions appear to depend on the choice of pseudohermitian structure, but we will show later that they are indeed independent of that choice and are therefore CR invariant conditions on $\hat{M}$. By Lemma \ref{covderivpi}, if the second fundamental form of $f$ vanishes, then $\hat{R}_\alpha^{\;\;a}=\hat{R}_{\alpha \;\; b}^{\;\;a\;\;b}$.\\
To prove Lemma \ref{2ndffthm}, we first assume the second fundamental form of F vanishes. Then all the forms in (\ref{omegaterms}) vanish. The vanishing of $\frac{i}{\hat{n}+2}\hat{R}_{\: \: \: \gamma}^a\theta^\gamma +iC^a\theta$ implies that $\hat{R}_{\: \: \: \gamma}^a = 0$. Combining this with the vanishing of $\hat{\omega}_\alpha^{\: \: \: a} +\frac{i}{\hat{n}+2}\hat{R}_\alpha^{\: \: \: a}\theta$ yields $\hat{\omega}_\alpha^{\: \: \: a} =0$ \\ 
The other direction of the proof is obvious.\\ QED.
\section{Chain Preserving CR Embeddings}
To prove Theorem \ref{gentargetchainthm} we will establish the following theorem.
\begin{thm}
Suppose $M$ and $\hat{M}$ are strictly pseudoconvex CR hypersurfaces and $f : M \hookrightarrow \hat{M}$ is a smooth CR embedding. Then $f$ sends chains on $M$ to chains on $\hat{M}$ if and only if
\begin{equation}
\left\{
\begin{array}{l}
\omega_{\alpha \;\;\;\beta}^{\;\;b} =0
\vspace{.2cm}
\\
\hat{S}_{a \;\;\alpha \bar{\beta}}^{\;\;a} = 0
\vspace{.2cm}
\\
\hat{S}_{\beta \;\;b}^{\;\;a\;\;b} = 0
\vspace{.2cm}
\\
\hat{V}_{\;\;\beta}^{a\;\;\beta}=0
\end{array}
\right.
\label{confcond}
\end{equation}
if and only if, with respect to any adapted coframe
\begin{equation}
\left\{
\begin{array}{l}
\omega_{\alpha \;\;\;\beta}^{\;\;b} =0
\vspace{.2cm}
\\
C_{\alpha \bar{\beta}}=0
\vspace{.2cm}
\\
\hat{R}_\alpha^{\;\;b} = 0
\vspace{.2cm}
\\
\hat{A}_{\;\;\; ; b}^{ab}-\hat{D}_{\;\;\; ; \bar{b}}^{a\bar{b}}=0
\end{array}
\right.
\label{hermcond}
\end{equation}
Moreover, $f$ sends chains on $M$ to chains on $\hat{M}$ if and only if the CR second fundamental form of $f$ vanishes and there is a lift of $f$ to a conformal isometry of associated Fefferman bundles with vanishing pseudo-Riemannian second fundamental form.
\label{gentargetthm}
\end{thm}

We consider a pair of coframes $\{ \theta, \theta^\alpha \}$, $\{ \hat{\theta}, \hat{\theta}^A \}$, adapted with respect to the embedding $f$, defined near a point $p \in M$ and $f(p) \in \hat{M}$ respectively. Fix an arbitrary transversal curve $\gamma \subset M$ through $p$ and consider the image $f(\gamma)$ through $f(p)$ in $\hat{M}$. Parameterize $\gamma$ so that $\gamma(0)=p$ and for simplification we identifty $M$ with its image in $\hat{M}$. There exists functions $a^\alpha (t)$, $\hat{a}^A(t)$ so that
\begin{equation}
\left\{
\begin{array}{cc}
\theta^\alpha = 2a^\alpha(t)\theta 
\vspace{.1cm}\\
\hat{\theta}^A = 2\hat{a}^A(t) \hat{\theta}
\end{array}
\right.
\label{thetais2a}
\end{equation}
along $\gamma$ and $f\circ \gamma$ respectively. By Theorem \ref{psxaltdef}, these curves are chains if, 
\begin{equation}
\left\{
\begin{array}{cc}
da^\alpha = 4ia^\alpha |a|^2 \theta -  a^\beta\phi_\beta^{\: \: \:\alpha}-\frac{1}{2} \phi^\alpha 
\vspace{.2cm}\\
d\hat{a}^A = 4i\hat{a}^A |\hat{a}|^2 \hat{\theta} -  \hat{a}^B\hat{\phi}_B^{\: \: \:A}-\frac{1}{2} \hat{\phi}^A
\end{array}
\right.
\label{adapteddiffeqn}
\end{equation}
on $M$ and $\hat{M}$ respectively. Since the coframes are adapted and the curve $\gamma$ is contained in $M$ we have $\hat{a}^\alpha = a^{\alpha}$ and $\hat{a}^b = 0$. Thus, the map $f$ sends chains on $M$ to chains on $\hat{M}$ if and only if 
\begin{equation}
\left\{
\begin{array}{l}
4ia^\alpha |a|^2 \theta -  a^\beta\phi_\beta^{\: \: \:\alpha}-\frac{1}{2} \phi^\alpha = 4i\hat{a}^\alpha |\hat{a}|^2 \hat{\theta} -  \hat{a}^\beta\hat{\phi}_\beta^{\: \: \:\alpha}-\frac{1}{2} \hat{\phi}^\alpha
\vspace{.2cm}\\
0 = 4i\hat{a}^b |\hat{a}|^2 \hat{\theta}  - \hat{a}^B\hat{\phi}_B^{\: \: \:b}-\frac{1}{2} \hat{\phi}^b
\end{array}
\right.
\label{chainstochains}
\end{equation}
for all curves through $p$ satisfying (\ref{thetais2a}). \\
Fixing a chain, using $a^\alpha=\hat{a}^\alpha$, $\hat{a}^b=0$, and $|a|^2=|\hat{a}|^2$, the first term in (\ref{chainstochains}) becomes
\begin{equation}
4ia^\alpha |a|^2 \theta -  a^\beta\phi_\beta^{\: \: \:\alpha}-\frac{1}{2} \phi^\alpha = 4ia^\alpha |a|^2 \theta -  a^\beta\hat{\phi}_\beta^{\: \: \:\alpha}-\frac{1}{2} \hat{\phi}^\alpha
\label{canceleqn}
\end{equation}
Recall equation (\ref{CMpullbacks}) which relates the pullback by $f$ of the Chern-Moser forms on $\hat{M}$:
\begin{align*}
\hat{\phi}_\beta^{\:\:\: \alpha} = \phi_\beta^{\:\:\:\alpha} + C_\beta^{\:\:\:\alpha}\theta, \: \: \: \hat{\phi}^\alpha = \phi^\alpha + C_\mu^{\:\:\: \alpha} \theta^\mu + F^\alpha \theta
\end{align*}
Plugging these into equation (\ref{canceleqn}), using (\ref{thetais2a}), and canceling, we then conclude that along the curve $\gamma$, (\ref{canceleqn}) holds if and only if
\begin{equation}
\begin{array}{rl}
a^\beta\phi_\beta^{\: \: \:\alpha}+\frac{1}{2} \phi^\alpha &= a^\beta( \phi_\beta^{\:\:\:\alpha} + C_\beta^{\:\:\:\alpha}\theta)+\frac{1}{2} (\phi^\alpha + C_\mu^{\:\:\: \alpha} \theta^\mu + F^\alpha \theta) \vspace{.2cm}\\\
\Leftrightarrow 0 &= a^\beta C_{\beta}^{\;\;\alpha}\theta+\frac{1}{2}C_\beta^{\;\;\alpha}\theta^\beta+\frac{1}{2}F^\alpha \theta \vspace{.2cm}\\
 &= (2a^\beta C_{\beta}^{\;\;\alpha}+\frac{1}{2}F^\alpha) \theta  \vspace{.2cm}\\
\Leftrightarrow 0 &= 2C_{\beta}^{\;\;\alpha}a^\beta+\frac{1}{2}F^\alpha \vspace{.2cm}\\
\end{array}
\label{polychaincond}
\end{equation}
where the last equivalence holds since the curve is transversal to the complex tangent space (ie $\theta (\dot{\gamma}) \neq 0$). This holds at all points in $\gamma$, in particular it holds when evaluated at the point $p \in \gamma \subset M$. As we vary the chain through $p$, the values of $a = (a^\alpha)$ sweep out $\mathbb{C}^n$ and thus equation (\ref{polychaincond}) holds as a polynomial of $n$ complex variables $a^\alpha$ with constant coefficients (evaluated at $p$). 
In \cite{EHZ}, it is shown that
\begin{equation}
F^\alpha=-\frac{2i}{n-1}(C_{\beta \;\; ; }^{\;\;\beta\;\;\alpha}-C_{\beta \;\; ;}^{\;\;\alpha\;\;\beta})
\end{equation}
Therefore $F^\alpha$ vanishes if $C_{\alpha \bar{\beta}}$ vanishes. We conclude that the first condition of (\ref{chainstochains}) holds for all chains if and only if $C_{\alpha \bar{\beta}}=0$. Observe that, by (\ref{iffeqn}), $f$ must then have a lift to a conformal isometry of the Fefferman bundles.\\
Again, evaluating along $\gamma$, using (\ref{CMWeb}) and using (\ref{thetais2a}), the second condition of (\ref{chainstochains}) holds if and only if
\[
\begin{array}{rl}
0 & = 4i\hat{a}^b |\hat{a}|^2 \hat{\theta}  - \hat{a}^B\hat{\phi}_B^{\: \: \:b}-\frac{1}{2} \hat{\phi}^b \vspace{.2cm}\\
& = -a^\beta\hat{\phi}_\beta^{\;\;\;b}-\frac{1}{2}\hat{\phi}^b \vspace{.2cm}\\
 & = -a^\beta(\hat{\omega}_\beta^{\;\;\;b}+\hat{D}_\beta^{\;\;\;b}\theta)-\frac{1}{2}(\hat{\tau}^b +\hat{D}_\beta^{\;\;\;b}\theta^\beta+\hat{E}^b\theta) \vspace{.2cm}\\
& = -a^\beta\hat{\omega}_{\beta\;\;\;\mu}^{\;\;\;b}\theta^\mu-a^\beta\hat{D}_\beta^{\;\;\;b}\theta-\frac{1}{2}\hat{D}_\beta^{\;\;\;b}\theta^\beta-\frac{1}{2}\hat{E}^b\theta \vspace{.2cm}\\
& = (-2{\omega}_{\beta\;\;\;\mu}^{\;\;\;b}a^\beta a^\mu\hat-2\hat{D}_\beta^{\;\;\;b}a^\beta-\frac{1}{2}\hat{E}^b)\theta \vspace{.2cm}\\
\Leftrightarrow 0 & = 2{\omega}_{\beta\;\;\;\mu}^{\;\;\;b}a^\beta a^\mu\hat-2\hat{D}_\beta^{\;\;\;b}a^\beta-\frac{1}{2}\hat{E}^b
\end{array}
\]
where again we have used $\theta (\dot{\gamma}) \neq 0$. Again, we vary the chain through $p$ so that the values of $a = (a^\alpha)$ sweep out $\mathbb{C}^n$ and thus the above equality holds as a degree 2 polynomial of $n$ complex variables $a^\alpha$ with constant coefficients (evaluated at $p$). Observe, by (\ref{DandE}), that
\begin{equation}
\begin{array}{rl}
\hat{D}_\beta^{\;\;\;b} &= \frac{i}{\hat{n}+2}\hat{R}_\beta^{\;\;\;b}-\frac{i}{2(\hat{n}+1)(\hat{n}+2)}\hat{R}g_\beta^{\;\;\;b}\vspace{.2cm}\\
& = \frac{i}{\hat{n}+2}\hat{R}_\beta^{\;\;\;b} \vspace{.2cm}\\
\hat{E}^b &= \frac{2i}{2\hat{n}+1}(\hat{A}_{\;\;\;\;;C}^{Cb}-\hat{D}_{\;\;\;\;;\bar{C}}^{b\bar{C}})
\end{array}
\label{normterms}
\end{equation}
So the second condition of (\ref{chainstochains}) holds for all chains through $p$ if and only if when evaluated at $p$ we have
\begin{equation}
\left\{
\begin{array}{c}
{\omega}_{\alpha\;\;\;\beta}^{\;\;\;b} = 0 \vspace{.2cm}\\
\hat{R}_\beta^{\;\;\;b} = 0 \vspace{.2cm}\\
\hat{A}_{\;\;\;\;;C}^{Cb}-\hat{D}_{\;\;\;\;;\bar{C}}^{b\bar{C}} = 0
\end{array}
\right.
\end{equation}
Since on $M$, $\hat{\tau}^b = \hat{A}_{\;\;\bar{\beta}}^{b}\theta^\alpha=0$, the covariant derivative is given by:
\begin{equation}
\begin{array}{rl}
\hat{\nabla}\hat{A}^{\beta b}&=d\hat{A}^{\beta b}+ \hat{A}^{C b}\hat{\omega}_{C}^{\;\;\;\beta}-\hat{A}^{\beta C}\hat{\omega}_{C}^{\;\;\;b}\\
&=\hat{A}^{a b}\hat{\omega}_{a}^{\;\;\;\beta}-\hat{A}^{\beta \mu}\hat{\omega}_{\mu}^{\;\;\;b}
\end{array}
\end{equation} 
Pulling back to $M$, if the second fundamental form of $f$ vanishes, then $\hat{\omega}_{a}^{\;\;\;\beta}=0$ and $\hat{\omega}_{\mu}^{\;\;\;b}=0$.Thus we conclude that $\hat{A}_{\;\;\;\; ; \alpha}^{\beta b}=0$, and thus $\hat{A}_{\;\;\;\;;\bar{C}}^{\bar{C}b}=\hat{A}_{\;\;\;\;;\bar{c}}^{\bar{c}b}$. Similarly, if $\hat{R}_\beta^{\;\;\;b} = 0$ on $M$ then $\hat{D}^{b \beta} = 0$ and by a similar computation we get $\hat{D}_{\;\;\;\;;\bar{C}}^{\bar{C}b}=\hat{D}_{\;\;\;\;;\bar{c}}^{\bar{c}b}$.\\
In summary, we have shown that the map $f$ sends chains through $p$ in $M$ to chains through $f(p)$ in $\hat{M}$ if and only if  (\ref{hermcond}) holds at $p$,
\[
\left\{
\begin{array}{l}
\omega_{\alpha \;\;\;\beta}^{\;\;b} =0
\vspace{.2cm}
\\
C_{\alpha \bar{\beta}}=0
\vspace{.2cm}
\\
\hat{R}_\alpha^{\;\;b} = 0
\vspace{.2cm}
\\
\hat{A}_{\;\;\; ; b}^{ab}-\hat{D}_{\;\;\; ; \bar{b}}^{a\bar{b}}=0
\end{array}
\right.
\]
We will now show that (\ref{hermcond}) is equivalent to (\ref{confcond}). First we will establish that
\begin{equation}
\left\{\begin{array}{c} \omega_{\alpha \;\;\;\beta}^{\;\;b} =0\vspace{.2cm}\\C_{\alpha \bar{\beta}}=0 \end{array}\right. 
\Leftrightarrow 
\left\{\begin{array}{c} \omega_{\alpha \;\;\;\beta}^{\;\;b} =0 \vspace{.2cm}\\\hat{S}_{a \;\;\alpha \bar{\beta}}^{\;\;a} = 0 \end{array}\right. 
\label{first2}
\end{equation}
Suppose $\omega_{\alpha \;\;\;\beta}^{\;\;b} =0$. In consideration of (\ref{Cs}) it is clear that $C_{\alpha \bar{\beta}}=0$ if and only if $\hat{S}_{a \;\;\alpha \bar{\beta}}^{\;\;a} - \frac{1}{2(n+1)}\hat{S}_{a \;\;\mu}^{\;\;a \;\;\mu}g_{\alpha \bar{\beta}}$. Contracting this gives $0=\hat{S}_{a \;\;\mu}^{\;\;a \;\; \mu} - \frac{1}{2(n+1)}\hat{S}_{a \;\;\mu}^{\;\;a \;\;\mu}g_{\alpha \bar{\beta}}=(1-\frac{n}{2(n+1)})\hat{S}_{a \;\;\mu}^{\;\;a \;\;\mu}=\frac{n+2}{2(n+1)}\hat{S}_{a \;\;\mu}^{\;\;a \;\;\mu}$. Thus we have shown that $C_{\alpha \bar{\beta}}=0$ if and only if $\hat{S}_{a \;\;\alpha \bar{\beta}}^{\;\;a}=0$ (when the CR second fundamental form of $f$ vanishes) and the equivalence in (\ref{first2}) is established. Thus, to finish proving Theorem \ref{gentargetthm}, we will show that
\begin{equation}
\left\{\begin{array}{c}\omega_{\alpha \;\;\;\beta}^{\;\;b} =0\vspace{.2cm}\\\hat{R}_\alpha^{\;\;b} = 0 \vspace{.2cm} \\ \hat{A}_{\;\;\; ; b}^{ab}-\hat{D}_{\;\;\; ; \bar{b}}^{a\bar{b}}=0  \end{array}\right. 
\Leftrightarrow
\left\{\begin{array}{c}\omega_{\alpha \;\;\;\beta}^{\;\;b} =0\vspace{.2cm}\\\hat{S}_{\beta \;\;b}^{\;\;a\;\;b} = 0 \vspace{.2cm} \\ \hat{V}_{\;\;\beta}^{a\;\;\beta}=0  \end{array}\right. 
\label{second2}
\end{equation}
The pseudohermitian curvature and pseudo-conformal curvature tensors are related by
\begin{equation}
\begin{array}{c}
\hat{S}_{A\bar{B}C\bar{D}}= \hat{R}_{A\bar{B}C\bar{D}}-
\frac{1}{\hat{n}+2}(\hat{R}_{A\bar{B}}g_{C\bar{D}}+\hat{R}_{C\bar{B}}g_{A\bar{D}}+\hat{R}_{A\bar{D}}g_{C\bar{B}}+\hat{R}_{C\bar{D}}g_{A\bar{B}})
\vspace{.2cm}\\
+\frac{\hat{R}}{(\hat{n}+1)(\hat{n}+2)}(g_{A\bar{B}}g_{C\bar{D}}+g_{A\bar{D}}g_{C\bar{B}})
\end{array}
\end{equation}
Using this we compute
\begin{equation}
\begin{array}{c}
\hat{S}_{\beta\;\;b}^{\;\;a\;\;b}= \hat{R}_{\beta\;\;b}^{\;\;a\;\;b}-
\frac{1}{\hat{n}+2}(\hat{R}_{\beta}^{\;\;a}g_{b}^{\;\;b}+\hat{R}_{\beta}^{\;\;b}g_{b}^{\;\;a}+\hat{R}_{b}^{\;\;a}g_{\beta}^{\;\;b}+\hat{R}_{b}^{\;\;b}g_{\beta}^{\;\;a})
\vspace{.2cm}\\
+\frac{\hat{R}}{(\hat{n}+1)(\hat{n}+2)}(g_{\beta}^{\;\;a}g_{b}^{\;\;b}+g_{\beta}^{\;\;b}g_{b}^{\;\;a})
\end{array}
\end{equation}
Since $g_{A\bar{B}}$ is the identity matrix, this gives
\begin{equation}
\begin{array}{c}
\hat{S}_{\beta\;\;b}^{\;\;a\;\;b}= \hat{R}_{\beta\;\;b}^{\;\;a\;\;b}-\frac{\hat{n}-n+1}{\hat{n}+2}\hat{R}_{\beta}^{\;\;a}
\end{array}
\end{equation}
By Lemma \ref{covlem}, when $\omega_{\alpha \;\;\;\beta}^{\;\;a} =0$, $\hat{R}_{\beta\;\;b}^{\;\;a\;\;b}=\hat{R}_{\beta}^{\;\;a}$, and thus
\begin{equation}
\begin{array}{c}
\hat{S}_{\beta\;\;b}^{\;\;a\;\;b}= \frac{n+1}{\hat{n}+2}\hat{R}_{\beta}^{\;\;a}
\end{array}
\end{equation}
In particular, $\hat{S}_{\beta\;\;b}^{\;\;a\;\;b}=0$ if and only if $\hat{R}_{\beta}^{\;\;a}$.
\vspace{.2cm}
\\
From (\ref{structureeqns})
\begin{equation}
\begin{array}{c}
d\hat{\phi}^A=\hat{\phi}\wedge\hat{\phi}^A+\hat{\phi}^B\wedge\hat{\phi}_B^{\;\;A}-\frac{1}{2}\hat{\psi}\wedge\hat{\omega}^A-\hat{V}_{B \; C}^{\;A}\hat{\omega}^B\wedge\hat{\omega}^C
\vspace{.1cm}\\
\;\;\;\;\;+\hat{V}_{\;\;B\bar{C}}^{A}\hat{\omega}^B\wedge\hat{\omega}^{\bar{C}}+\hat{P}_B^{\;\;A}\hat{\omega}^B\wedge\hat{\omega}+\hat{Q}_{\bar{C}}^{\;\;A}\hat{\omega}^{\bar{C}}\wedge\hat{\omega}
\end{array}
\end{equation}
Since $\hat{\phi}=0$, $\hat{\omega}^A=\hat{\theta}^A$, $\hat{A}_{\;\;\bar{\mu}}^b=0$ setting $A=a$, using the analog of (\ref{CMWeb}) on $\hat{M}$, and pulling back to $M$, this becomes
\begin{equation}
\begin{array}{rl}
d\hat{\phi}^a&=\hat{\phi}^B\wedge\hat{\phi}_B^{\;\;a}-\hat{V}_{\alpha \; \gamma}^{\;a}\theta^\alpha\wedge\theta^\gamma+\hat{V}_{\;\;\alpha\bar{\mu}}^{a}\theta^\alpha\wedge\theta^{\bar{\mu}}+\hat{P}_\alpha^{\;\;a}\theta^\alpha\wedge\theta+\hat{Q}_{\bar{\mu}}^{\;\;a}\theta^{\bar{\mu}}\wedge\theta 
\vspace{.2cm}\\
&=(\hat{A}_{\;\;\bar{\mu}}^B\theta^{\bar{\mu}} + \hat{D}_{\alpha}^{\: \: \: B}\theta^\alpha + \hat{E}^B\theta )\wedge(\hat{\omega}_{B}^{\: \: \: a} + D_{\beta}^{\: \: \: \alpha}\theta)-\hat{V}_{\alpha \; \gamma}^{\;a}\theta^\alpha\wedge\theta^\gamma
\vspace{.1cm}\\
&\;\;\;\;\;\;+\hat{V}_{\;\;\alpha\bar{\mu}}^{a}\theta^\alpha\wedge\theta^{\bar{\mu}}+\hat{P}_\alpha^{\;\;a}\theta^\alpha\wedge\theta+\hat{Q}_{\bar{\mu}}^{\;\;a}\theta^{\bar{\mu}}\wedge\theta
\vspace{.2cm}\\
&= (\hat{V}_{\;\;\;\alpha\bar{\mu}}^{a}-\hat{A}_{\;\;\bar{\mu}}^{\beta}\omega_{\beta \;\; \alpha}^{\;\;a})\theta^\alpha \wedge \theta^{\bar{\mu}}-\hat{V}_{\alpha\;\;\gamma}^{\;\;a}\theta^\alpha\wedge\theta^\gamma+(\hat{D}_\alpha^{\;\;B}\theta^\alpha+\hat{E}^B\theta)\wedge\hat{\omega}_B^{\;\;a}
\vspace{.1cm}\\
&\;\;\;\;\;+(\hat{A}_{\;\;\bar{\mu}}^{\beta}\hat{D}_\beta^{\;\;a}+\hat{Q}_{\bar{\mu}}^{\;\;a})\theta^{\bar{\mu}}\wedge\theta+(\hat{D}_\alpha^{\;\;\beta}\hat{D}_\beta^{\;\;a}+\hat{P}_\alpha^{\;\;a})\theta^\alpha\wedge\theta
\end{array}
\label{dtheta1}
\end{equation}
On the other hand, $\hat{\phi}^a=\hat{\tau}^a+\hat{D}_C^{\;\;a}\theta^C+\hat{E}^a\theta=\hat{D}_\alpha^{\;\;a}\theta^\alpha+\hat{E}^a\theta$, so
\begin{equation}
\begin{array}{rl}
d\hat{\phi}^a&=d\hat{D}_C^{\;\;a}\wedge\theta^C+\hat{D}_C^{\;\;a}d\theta^C+d\hat{E}^a\wedge\theta+\hat{E}^ad\theta
\vspace{.2cm}\\
& = d\hat{D}_\alpha^{\;\;a}\wedge\theta^\alpha+\hat{D}_C^{\;\;a}(\theta^\alpha\wedge\omega_\alpha^{\;\;C}+\theta\wedge\tau^C)+d\hat{E}^a\wedge\theta+i\hat{E}^ag_{\alpha\bar{\beta}}\theta^\alpha\wedge\theta^{\bar{\beta}}
\vspace{.2cm}\\
& = (d\hat{D}_\alpha^{\;\;a}-\hat{D}_C^{\;\;a}\hat{\omega}_\alpha^{\;\;C})\wedge\theta^\alpha+i\hat{E}^ag_{\alpha\bar{\beta}}\theta^\alpha\wedge\theta^{\bar{\beta}}+(d\hat{E}^a-\hat{D}_\alpha^{\;\;a} A_{\;\;\bar{\mu}}^\alpha\theta^{\bar{\mu}})\wedge\theta
\end{array}
\label{dtheta2}
\end{equation}
Combining (\ref{dtheta1}) and (\ref{dtheta2}), we conclude
\begin{equation}
\begin{array}{l} i\hat{E}^ag_{\alpha\bar{\beta}}\theta^\alpha\wedge\theta^{\bar{\beta}}+\hat{\nabla}\hat{D}_\alpha^{\;\;a}\wedge\theta^\alpha+(\hat{\nabla}\hat{E}^a-\hat{D}_\alpha^{\;\;a} A_{\;\;\bar{\mu}}^\alpha\theta^{\bar{\mu}})\wedge\theta
\vspace{.2cm}\\
\;\;=(\hat{V}_{\;\;\;\alpha\bar{\mu}}^{a}-\hat{A}_{\;\;\bar{\mu}}^{\beta}\omega_{\beta \;\; \alpha}^{\;\;a})\theta^\alpha \wedge \theta^{\bar{\mu}}-\hat{V}_{\alpha\;\;\gamma}^{\;\;a}\theta^\alpha\wedge\theta^\gamma
\vspace{.1cm}\\
\;\;\;\;\;\;\;\;+(\hat{A}_{\;\;\bar{\mu}}^{\beta}\hat{D}_\beta^{\;\;a}+\hat{Q}_{\bar{\mu}}^{\;\;a})\theta^{\bar{\mu}}\wedge\theta+(\hat{D}_\alpha^{\;\;\beta}\hat{D}_\beta^{\;\;a}+\hat{P}_\alpha^{\;\;a})\theta^\alpha\wedge\theta
\end{array}
\label{generalcovequality}
\end{equation}
If $\omega_{\alpha \;\; \beta}^{\;\;a}=0$ and $\hat{R}_\alpha^{\;\;a}=0$ (so then $\hat{D}_\alpha^{\;\;\;a}=0$) then we have $\hat{\nabla}\hat{D}_\alpha^{\;\;a}=0$. Comparing similar terms in (\ref{generalcovequality}) we conclude that 
\begin{equation}
i\hat{E}^ag_{\alpha\bar{\beta}}=\hat{V}_{\;\;\;\alpha\bar{\beta}}^a
\end{equation}
contracting both sides we get
\begin{equation}
\hat{E}^a=-\frac{i}{n}\hat{V}_{\;\;\beta}^{a\;\;\beta}
\end{equation}
and thus
\begin{equation}
\hat{E}^a = 0 \text{ if and only if } \hat{V}_{\;\;\beta}^{a\;\;\beta} = 0
\end{equation}
By the comments following (\ref{normterms}), since $\omega_{\alpha \;\; \beta}^{\;\;a}=0$ and $\hat{R}_\alpha^{\;\;a}=0$, 
\begin{equation}
\hat{E}^a = \frac{2i}{2\hat{n}+1}(\hat{A}_{\;\;\;\;;b}^{ab}-\hat{D}_{\;\;\;\;;\bar{b}}^{a\bar{b}})
\end{equation}
Note that, by the trace condition $\hat{V}_{\;\;C}^{a\;\;C}=0$, $\hat{V}_{\;\;\beta}^{a\;\;\beta}=-\hat{V}_{\;\;b}^{a\;\;b}$.
We have now shown that (\ref{hermcond}) is equivalent to (\ref{confcond}). To finish Theorem \ref{gentargetthm}, we will show that the conditions in (\ref{hermcond}) hold if and only if there is a lift of $f$ to a conformal isometry of the Fefferman bundles with vanishing pseudo-Riemannian second fundamental form. \\
We know from (\ref{iffeqn}) that $C_{\alpha\bar{\beta}}=0$ if and only if there is a lift of $f$ to a conformal isometry of the Fefferman bundles. Comparing Lemma \ref{2ndffthm} and (\ref{hermcond}) we see it suffices to show that when $\omega_{\alpha \;\; \beta}^{\;\;a}=0$ and $\hat{R}_\alpha^{\;\;a}=0$, then
\begin{equation}
\hat{A}_{\bar{b}\;\; ;}^{\;\;\;a\;\bar{b}} -\frac{i}{2(\hat{n}+1)}\hat{R}_{\: ;}^{\: \: a} = 0 \text{ if and only if } \hat{A}_{\;\;\;\;;b}^{ab}-\hat{D}_{\;\;\;\;;\bar{b}}^{a\bar{b}}=0
\label{finalequality}
\end{equation}
By one of the Bianchi identities given in \cite{Lee88}
\begin{equation}
\begin{array}{rl}
\hat{R}_{\;\;C;}^{a\;\;\;C}&=\hat{R}_{\: ;}^{\: \: a}-i(\hat{n}-1)\hat{A}_{\bar{C}\;\;;}^{\;\;\;a\;\bar{C}}
\end{array}
\end{equation}
using the assumptions that  $\omega_{\alpha \;\; \beta}^{\;\;a}=0$ and $\hat{R}_\alpha^{\;\;a}=0$ we have $\hat{R}_{\;\;C;}^{a\;\;\;C} =\hat{R}_{\;\;b;}^{a\;\;\;b}$ and $\hat{A}_{\bar{C}\;\;;}^{\;\;\;a\;\bar{C}}=\hat{A}_{\bar{b}\;\;;}^{\;\;\;a\;\bar{b}}$, thus
\begin{equation}
\hat{R}_{\;\;b;}^{a\;\;\;b}=\hat{R}_{\: ;}^{\: \: a}+i(\hat{n}-1)\hat{A}_{\bar{b}\;\;;}^{\;\;\;a\;\bar{b}}
\end{equation}
Notice that, by raising an index, we have introduced a conjugate. This becomes
\begin{equation}
\begin{array}{rl}
\hat{R}_{\;\;\;\;;\bar{b}}^{a\bar{b}}= \hat{R}_{\: ;}^{\: \: a}+i(\hat{n}-1)\hat{A}_{\bar{b}\;\;;}^{\;\;\;a\;\bar{b}}
\end{array}
\label{bianchiscalar}
\end{equation}
By definition and using (\ref{bianchiscalar}):
\begin{equation}
\begin{array}{rl}
\hat{D}_{\;\;\;\;;\bar{b}}^{a\bar{b}}	&	= \frac{i}{\hat{n}+1}(\hat{R}^{a\bar{b}}-\frac{1}{2(\hat{n}+1)}\hat{R}g^{a\bar{b}})_{;\bar{b}}
\vspace{.2cm}\\
& = \frac{i}{\hat{n}+1}(\hat{R}_{\;\;\;\;;\bar{b}}^{a\bar{b}}-\frac{1}{2(\hat{n}+1)}\hat{R}_{;}^{\;a})
\vspace{.2cm}\\
& = \frac{i}{\hat{n}+1}(\hat{R}_{\: ;}^{\: \: a}+i(\hat{n}-1)\hat{A}_{\bar{b}\;\;;}^{\;\;\;a\;\bar{b}}-\frac{1}{2(\hat{n}+1)}\hat{R}_{;}^{\;a})
\vspace{.2cm}\\
& = \frac{i}{\hat{n}+1}(\frac{2\hat{n}+1}{2(\hat{n}+1)}\hat{R}_{\: ;}^{\: \: a}+i(\hat{n}-1)\hat{A}_{\bar{b}\;\;;}^{\;\;\;a\;\bar{b}})
\end{array}
\end{equation}
Thus
\begin{equation}
\begin{array}{rl}
\hat{A}_{\;\;\;\;;b}^{ab}-\hat{D}_{\;\;\;\;;\bar{b}}^{a\bar{b}} & = \hat{A}_{\;\;\;\;;b}^{ab}-\frac{i}{\hat{n}+1}(\frac{2\hat{n}+1}{2(\hat{n}+1)}\hat{R}_{\: ;}^{\: \: a}+i(\hat{n}-1)\hat{A}_{\bar{b}\;\;;}^{\;\;\;a\;\bar{b}})
\vspace{.2cm}\\
& = \frac{2\hat{n}+1}{\hat{n}+2}(\hat{A}_{\bar{b}\;\;;}^{\;\;a\;\bar{b}}-\frac{i}{2(\hat{n}+1)}\hat{R}_{\: ;}^{\: \: a})
\end{array}
\end{equation}
We have shown that the statements in (\ref{finalequality}) are equivalent and thus, $f$ sends chains on $M$ to chains on $\hat{M}$ if and only if there is a lift of $f$ to a conformal isometry between the Fefferman bundles with vanishing pseudo-Riemannian second fundamental form. 
\section{Mappings Into Spheres}
Now we assume $\hat{M}$ is (locally) the sphere $\mathbb{S}^{2\hat{n}+1} \subset \mathbb{C}^{\hat{n}+1}$ and thus its CR curvature tensor vanishes identically. We will prove
\begin{thm}
Suppose $f : M \hookrightarrow \mathbb{S}^{2\hat{n}+1}$ is a local smooth CR embedding of a strictly pseudoconvex smooth hypersurface $M \subset \mathbb{C}^{n+1}$ into the sphere $\mathbb{S}^{2\hat{n}+1}\subset\mathbb{C}^{\hat{n}+1}$. Let $C$ and $\hat{C}$ be the Fefferman bundles associated to $M$ and $\mathbb{S}^{2\hat{n}+1}$ respectively. The following conditions are equivalent.
\begin{enumerate}
\item $f$ sends chains on $M$ to chains on $\hat{M}$.
\vspace{.1cm}
\item There is a lift of $f$ to a conformal isometry between the Fefferman bundles $C$ and $\hat{C}$.
\vspace{.1cm}
\item The CR second fundamental form of f vanishes. 
\vspace{.1cm}
\item There exists a local CR diffeomorphism $\phi$ from the sphere $\mathbb{S}^{2n+1}$ to $M$ and an automorphism of the target sphere $A \in Aut(\mathbb{S}^{2\hat{n}+1})$ such that the composition $A \circ f \circ \phi$ : $ \mathbb{S}^{2n+1} \rightarrow \mathbb{S}^{2\hat{n}+1}$ is the linear embedding.
\end{enumerate}
\label{fefthm2}
\end{thm}
$1 \Rightarrow 2$ and $1 \Rightarrow 3$ in Theorem \ref{fefthm2} holds by the previous section. To prove $2 \Leftrightarrow 3$ we observe by (\ref{techeqn}) that the map $f$ may be lifted to a conformal isometry between Fefferman metrics if and only if
\begin{equation}
\omega_{\mu \: \: \: \alpha}^{\: \: \: a}\omega_{\: \: \: a \bar{\beta} }^\mu = 0
\label{sphereeqn}
\end{equation}
\begin{lem}
Condition (\ref{sphereeqn}) is satisfied if and only if $\omega_{\alpha \: \: \: \beta}^{\:\:\:a} = 0$
\end{lem}
The 'if' portion of the lemma is obvious. Now let us assume (\ref{sphereeqn}) holds. Contracting $\alpha$ and $\bar{\beta}$ we see
\begin{equation}
0=\omega_{\mu \: \: \: \nu}^{\: \: \: a}\omega_{\: \: \: a  }^{\mu\;\;\;\nu} = \displaystyle \sum_{\mu,\nu,a}|\omega_{\mu \: \: \: \nu}^{\: \: \: a}|^2
\end{equation}
QED 
\vspace{.5cm}
\\
Now to prove $3 \Rightarrow 4$ in Theorem \ref{fefthmsphere} we first recall that in \cite{EHZ} the following pseudo-conformal Gauss equation is established for all $p\in M$:
\begin{equation}
[ \hat{S}(X,Y,Z,V) ] = S(X,Y,Z,V)+ [ \langle \Pi(X,Z),\Pi(Y,V) \rangle ], \: \: \: \: X,Y,Z,V \in T_p^{(1,0)}M
\end{equation}
where $[\hat{S}]$ denotes the traceless component of the pseudo-conformal curvature tensor $\hat{S}$ on $\hat{M}$. Since both $\hat{S}$ and $\Pi$ vanish identically clearly so must $S$, that is to say that $M$ is CR-flat and thus locally equivalent to the sphere. In \cite{EHZ} the following rigidity result is proved:
\begin{thm}
Let $f : M \hookrightarrow \mathbb{S}^{2\hat{n}+1}$ be a smooth CR-immersion and s be the degeneracy of $f$. If $\hat{n}-n-s \leq \frac{n}{2}$, then any other such CR-immersion $\tilde{f}$ is related to $f$ by $\tilde{f}=A \circ f$, where $A$ is a CR-automorphism of the ambient sphere. 
\end{thm}
Since the CR second fundamental form of $f$ vanishes, its degeneracy is $s=\hat{n}-n$. Let $\phi$ : $\mathbb{S}^{2n+1} \rightarrow M$ be a local CR diffeomorphism guaranteed now since $M$ is CR flat. The degeneracy of $f \circ \phi$ is the same as that of $f$ and the CR second fundamental form of $f \circ \phi$ still vanishes, thus $f \circ \phi$ is equivalent to the trivial map, ie there is an automorphism $A \in Aut(\mathbb{S}^{2\hat{n}+1})$ such that the composition $A \circ f \circ \phi$ : $ \mathbb{S}^{2n+1} \rightarrow \mathbb{S}^{2\hat{n}+1}$ is the linear embedding.\\
Now we finish the proof of Theorem \ref{fefthmsphere} by concluding $4 \Rightarrow 1$. Suppose maps $\phi:\mathbb{S}^{2n+1} \rightarrow M$ and $A \in Aut(\mathbb{S}^{2\hat{n}+1})$ exist as in the theorem so that $A\circ f \circ \phi :\mathbb{S}^{2n+1} \rightarrow \mathbb{S}^{2\hat{n}+1}$ is linear. It is known that the chains in a sphere are exactly the great circles and thus a linear embedding between spheres preserves chains. CR diffeomorphisms always preserve chains and thus since $\phi$, $A$ and $A\circ f \circ \phi$ all preserve chains, $f$ must as well. Equivalently, it is easy to see that the conditions in (\ref{confcond}) are satisfied. 
\\
We have now established Theorem \ref{fefthmsphere}.
\newpage
\bibliographystyle{unsrt}
\bibliography{mybib}
\vspace{.5cm}
A. Minor, aminor@ucsd.edu, Department of Mathematics, University of California San Diego, La Jolla, CA, 92093, USA.
\end{document}